\documentclass[11pt]{amsart}
\usepackage[T1]{fontenc}
\usepackage{lmodern}
\usepackage[utf8]{inputenc}
\usepackage{amsmath,amssymb,amsthm,mathrsfs}
\usepackage{cite}
\usepackage[colorlinks=true,linkcolor=blue,citecolor=blue,urlcolor=blue]{hyperref}
\usepackage[nameinlink,noabbrev]{cleveref}
\numberwithin{equation}{section}
\theoremstyle{plain}
\newtheorem{theorem}{Theorem}[section]
\newtheorem{corollary}[theorem]{Corollary}
\newtheorem{lemma}[theorem]{Lemma}
\newtheorem{proposition}[theorem]{Proposition}
\theoremstyle{definition}

\newtheorem{problem}[theorem]{Problem}
\theoremstyle{remark}
\newtheorem{remark}[theorem]{Remark}
\newcommand{\C}{\mathbb C}
\newcommand{\R}{\mathbb R}
\newcommand{\D}{\Delta}
\newcommand{\T}{\mathcal T}
\newcommand{\OO}{\mathcal O}
\newcommand{\rank}{\operatorname{rank}}
\newcommand{\Span}{\operatorname{span}}
\newcommand{\Ann}{\operatorname{Ann}}
\newcommand{\im}{\operatorname{im}}
\newcommand{\id}{\operatorname{id}}
\newcommand{\dd}{\mathrm d}
\newcommand{\Levi}{\mathcal L}
\newcommand{\Rea}{\operatorname{Re}}
\newcommand{\Ima}{\operatorname{Im}}

\Crefname{problem}{Problem}{Problems}
\Crefname{remark}{Remark}{Remarks}
\date{}
\title{Local convexity obstructions for bounded realizations of Teichm\"uller space}
\author{Jinsong Liu\textsuperscript{1} $\&$ Xingsi Pu\textsuperscript{2} $\&$ Lang Wang\textsuperscript{3} }
\address{$1.$ State Key Laboratory of Mathematical Sciences, AMSS, Chinese Academy of Sciences, Beijing, 100190, China; School of Mathematical Sciences, University of Chinese Academy of Sciences, Beijing, 100049, China}
\address{$2.$ Mathematical Science Research Center, Chongqing University of Technology, Chongqing, 400054, China}
\address{$3.$ School of Mathematical Sciences, Guizhou Normal University, Guiyang, 550025, P.R. China.}
\email{liujsong@math.ac.cn,\:puxs@cqut.edu.cn,\:wanglang2020@amss.ac.cn}


\begin{document}
\begin{abstract}
Let $D\subset\C^n$, $n\ge2$, be a bounded domain that is biholomorphic to a
finite-dimensional Teichm\"uller space. We prove that $\partial D$ has no
locally convex boundary patch containing a complex $(n-1)$-dimensional
manifold through each of its points. Note that in complex dimension two this
excludes every locally convex boundary point, without any boundary
regularity assumption. In every dimension, we also exclude locally
convex $C^2$ boundary patches.
\end{abstract}
\maketitle

\section{Introduction}\label{sec:int}

Let $S_{g,m}$ be an oriented surface of genus $g$ with $m$ punctures,
where $2g-2+m>0$, and let $\T_{g,m}$ be its Teichm\"uller space. Throughout this paper,
\[
 n=\dim_{\C}\T_{g,m}=3g-3+m\ge2,
 \quad \Psi:\T_{g,m}\longrightarrow D
\]
is a biholomorphism, where $D\subset\C^n$ is bounded. Bers' embedding realizes
Teichm\"uller space as a bounded domain in $\C^n$, and its complex-analytic and metric
background is described in \cite{Bers,GL}.
For closed surfaces of genus at least $2$, Kim \cite{Kim} proved that the image of the Bers embedding is not convex.
In this paper, we will study local convexity obstructions of  arbitrary bounded realizations of finite-dimensional Teichm\"uller spaces.

The global convexity obstruction arises as follows. Markovic \cite{Mark}  first proved that the Carath\'eodory and Kobayashi metrics on
the Teichm\"uller space of a closed surface of genus at least two do not
coincide. Gekhtman and Markovic \cite{GM} extended this
noncoincidence to all finite-type Teichm\"uller spaces of complex dimension at
least two. Combined with Lempert's theorem \cite{Lem}, this rules out bounded
convex realizations.

On the other hand, Gupta and Seshadri \cite{GS} excluded locally strictly
convex boundary points for realizations of closed-surface Teichm\"uller
spaces.  Liu and Wang \cite{LW} further excluded
locally log-type convex boundary points.
In our earlier work \cite{PW}, we excluded bounded realizations with a good corner point
and also proved a local peak-point obstruction.
That is, no such bounded realization can
be locally convex at a boundary point admitting a local holomorphic
peak function.

For a bounded convex domain, a boundary point is a holomorphic peak point
if and only if no nonconstant analytic disc in the boundary passes
through that point \cite[Theorem 6]{NP}.
Applied to a convex local intersection, it shows that every point of a
locally convex patch of a bounded realization must lie on such a disc.
This leaves the following question.

\begin{problem}
Can a bounded realization of $\T_{g,m}$ have a locally convex
boundary patch such that every point of the patch lies on a nonconstant
analytic disc contained locally in the boundary?
\end{problem}

A negative answer would exclude all locally convex boundary points.
We first treat the case in which every point lies on a complex
$(n-1)$-dimensional manifold contained locally in the boundary.

\begin{theorem}\label{thm:nonsmooth}
There is no nonempty locally convex boundary patch $M$ of $D$ such that,
for every $q\in M$, there is a holomorphic map
\[
 f_q:\mathbb{B}^{n-1}\longrightarrow M,\quad
 f_q(0)=q,\quad \rank_{\C}\dd f_q(0)=n-1.
\]
Here $\mathbb{B}^k$ denotes the unit ball in $\C^k$. The rank hypothesis means
that each map parametrizes an embedded complex manifold after restriction
to a sufficiently small ball about its center.
\end{theorem}

Combining this result with the local peak-point obstruction gives a
complete answer in complex dimension $2$.
\begin{corollary}\label{cor:disc-free}
A bounded realization $D\subset\C^2$ of a two-dimensional
Teichm\"uller space has no locally convex boundary point.
\end{corollary}

We also exclude $C^2$ patches in every dimension. Note that here the boundary may
have positive Levi rank and nontrivial complex flat directions.

\begin{theorem}\label{cl:main}
The boundary of a bounded realization $D \subset \C^n$ as above has no nonempty locally
convex $C^2$ boundary patch.
\end{theorem}

The paper is organized as follows. In \Cref{sec2}, we review disc rigidity and the relevant results in convex geometry. In \Cref{sec:nonsmooth}, we prove \Cref{thm:nonsmooth} and \Cref{cor:disc-free}. In \Cref{sec:analytic-tools}, we develop the analytic tools and provide the direct affine-fiber construction. In \Cref{sec:conormal}, we establish the conormal reduction for the remaining case. Finally, \Cref{sec:family} completes the proof of \Cref{cl:main}.

\section{Preliminaries}\label{sec2}
\subsection{Notation}
\begin{enumerate}
    \item The Euclidean complex structure is $Jv=iv$, and $\langle u,v\rangle=\Rea\sum_j\overline{u_j}v_j$ is the real inner product.

    \item Write
    $$
    \D=\{\lambda\in\C:|\lambda|<1\}, \:
    \D_\epsilon=\{s\in\C:|s|<\epsilon\},
    $$
    and $\mathbb{B}_\epsilon^k=\{s\in\C^k:|s|<\epsilon\}$.
    We use $\OO$ for holomorphic maps or functions.

    \item A boundary analytic disc through $p$ is a nonconstant
    holomorphic map $f:\D\to\partial D$ with $f(0)=p$.
    Unless a rank condition is stated, its derivative at $0$ is allowed
    to vanish. All boundary discs considered locally can be restricted
    to a smaller disc.
\end{enumerate}

Recall that a boundary point $p\in\partial D$ is \emph{locally convex} if
$D\cap U$ is convex for some Euclidean ball $U$ centered at $p$.
We call $M=\partial D\cap U$ the true boundary patch, to distinguish
it from the boundary introduced by $\partial U$.
Moreover, a \emph{locally convex $C^2$ boundary patch} means that $M$
is an embedded $C^2$ real hypersurface.
We always work subsequently in a
smaller neighborhood compactly contained in $U$.

A \emph{local holomorphic peak function} at $p$ is a function
$h\in\OO(D\cap V)\cap C^0(\overline D\cap V)$ such that $h(p)=1$ and $|h|<1$ on
$(\overline D\cap V)\setminus\{p\}$, where $V$ is a neighborhood of $p$.
\begin{remark}
Note that the proof of \cite[Theorem 1.15]{PW} applies equally to finite-type punctured surfaces and to local peak functions in the sense just defined.
We use this form of the obstruction throughout.
\end{remark}

\subsection{Analytic discs and rigidity theorem}\label{sec:rigidity}
We recall the following peak-point criterion in the form given in \cite[Theorem 6]{NP}. The nontrivial implication, namely that the absence of affine complex discs implies the existence of a peak function, is due to Sibony \cite{S87}.

\begin{lemma}\label{lem:peak-criterion}
Let $\Omega\subset\C^n$ be a bounded convex domain and $p\in\partial\Omega$.
The following are equivalent:
\begin{enumerate}
 \item $p$ is a peak point for
 $A(\Omega)=\OO(\Omega)\cap C^0(\overline \Omega)$;
 \item there is no nonconstant holomorphic map $f:\D\to\partial \Omega$ with
 $f(0)=p$;
 \item there are no $v\ne0$ and $\epsilon>0$ for which
 $p+\lambda v\in\partial\Omega$ whenever $|\lambda|<\epsilon$.
\end{enumerate}
In particular, at a locally convex point of any domain, the absence of
a boundary analytic disc germ implies the existence of a local
holomorphic peak function.
\end{lemma}

\begin{lemma}\label{lem:straight}
Let $\Omega\subset\C^n$ be a convex domain and $p\in\partial \Omega$. Suppose
that $f:\mathbb{B}^{n-1}\to\partial \Omega$ is a holomorphic mapping with $f(0)=p$ and
$\rank_{\C}\dd f(0)=n-1$.
Then $\partial \Omega$ contains a relative neighborhood of $p$ in an affine complex hyperplane $H_p$. The direction of $H_p$ is unique, and every real supporting functional at $p$ annihilates this complex direction.
If $n=2$, the same conclusion holds for any nonconstant holomorphic disc through $p$, even if its derivative vanishes at the center.
\end{lemma}
\begin{proof}
Choose a nonzero complex-linear functional $\ell$ such that
$\Rea\ell(z-p)\le 0$ on $\overline\Omega$. The harmonic function
$\Rea\ell(f-p)$ attains an interior maximum, hence is constant.
Holomorphicity and the value at $0$ then yield $\ell(f-p)=0$.
The affine subspace $p+\ker\ell$ has complex dimension $n-1$, so the
holomorphic inverse function theorem implies that the image is relatively
open near $p$.

Every other supporting functional also annihilates the image, and hence
annihilates its $(n-1)$-dimensional complex direction. This proves
uniqueness. For $n=2$, the map takes values in a complex line, and the open
mapping theorem supplies a neighborhood of $p$ in that line.
\end{proof}

Recall that the Teichm\"uller metric $d_{\T}$ is a complete geodesic metric and agrees
with the Kobayashi distance \cite{Roy,GL}. Finite-dimensional Teichm\"uller
space is locally compact. The metric Hopf--Rinow theorem therefore shows
that $d_{\T}$ is proper.
We will use the following consequence of the rigidity theorem for holomorphic families \cite[Theorem 1]{Tan}, as well as the positive-measure unbounded-ray form in \cite[Theorem 2.2]{Ant}.

\begin{proposition}\label{prop:forbidden}
	There is no holomorphic map
	$F:\D_\epsilon\times\D\to D$
	whose slices are not all identical and for which
	$F(0,\cdot)$ extends continuously to a nonempty open arc
	$I\subset\partial\D$ with $F(0,I)\subset\partial D$.
\end{proposition}

\begin{proof}
  Suppose such a family exists. Denote
	\[
	G(\tau,\lambda)=\Psi^{-1}\bigl(F(\epsilon\tau,\lambda)\bigr),
	\quad (\tau,\lambda)\in\D\times\D,
	\]
	and write $g_\tau(\lambda)=G(\tau,\lambda)$.
	For $\xi\in I$, the points $F(0,t\xi)$ converge to a point of
	$\partial D$ as $t\rightarrow1$.
	Consequently $g_0(t\xi)$ leaves every compact subset of $\T_{g,m}$.
	Otherwise a subsequence in a compact set $L\subset\T_{g,m}$ would
	place the corresponding points in the compact set
	$\Psi(L)\subset D$.
	Then the properness of Teichm\"uller distance gives that
	\[
	d_{\T}\bigl(g_0(0),g_0(t\xi)\bigr)\longrightarrow\infty.
	\]
	Thus the central disc has unbounded radial image on a set $I$ of
	positive arc-length measure.
The positive-measure rigidity statement in \cite[Theorem 2.2]{Ant} therefore implies $g_\tau=g_0$ for all $\tau$.
	Since $\Psi$ is injective, all slices of $F$ are identical, which gives a
	contradiction.
\end{proof}

\section{Proof of \Cref{thm:nonsmooth}}\label{sec:nonsmooth}

\begin{proposition}
\label{prop:nonsmooth-product}
Let $\Omega\subset\mathbb C^n$ be a convex domain, and let
$M\subset\partial\Omega$ be nonempty and relatively open. Suppose that every
point of $M$ lies in the relative interior of an affine complex
$(n-1)$-dimensional plaque contained in $\partial\Omega$.
Then, after a complex-affine change of coordinates
$(w,z')\in\mathbb C\times\mathbb C^{n-1}$, there exist a ball
$B'\subset\mathbb C^{n-1}$ centered at $0$ and an open rectangle
$W=(-h,h)+i(-\eta,\eta)$ with $h,\eta>0$ such that
\[
 \partial\Omega\cap(W\times B')=\Gamma\times B'\subset M,
 \quad
 \Omega\cap(W\times B')=G\times B',
\]
where $\Gamma$ is a Lipschitz arc and $G\subset W$ is a planar domain
with $\partial G\cap W=\Gamma$. More precisely, for a convex Lipschitz
function $g:(-\eta,\eta)\to(-h,h)$ with $g(0)=0$, one can take
\[
 \Gamma=\{-g(t)+it:|t|<\eta\},
 \quad
 G=\{x+it\in W:x<-g(t)\}.
\]
\end{proposition}

\begin{proof}
We first find a uniform plaque radius on a nonempty open patch.
Choose a compact boundary graph patch $K\subset M$ such that
\[
 K^\circ:=\operatorname{int}_{\partial\Omega}K\ne\varnothing,
 \quad K=\overline{K^\circ}.
\]
For each positive integer $j$, let
\[
 K_j=\left\{q\in K:
 \begin{array}{l}
 \text{there is a complex linear hyperplane }H\subset\mathbb C^n
 \\[-2pt]
 \text{such that }q+
       \bigl(H\cap\overline{B(0,1/j)}\bigr)\subset\partial\Omega
 \end{array}\right\}.
\]
The plaque hypothesis gives $K=\bigcup_{j\ge1}K_j$.
Each $K_j$ is closed in $K$.

Indeed, if $q_\nu\in K_j$ tends to $q$,
choose corresponding hyperplanes $H_\nu$ and pass to a subsequence
converging in the complex Grassmannian to a hyperplane $H$.
For $v\in H\cap\overline{B(0,1/j)}$, the orthogonal projections
$v_\nu=P_{H_\nu}v$ satisfy $v_\nu\to v$ and $|v_\nu|\le1/j$.
Since $q_\nu+v_\nu\in\partial\Omega$ and $\partial\Omega$ is closed,
we have $q+v\in\partial\Omega$.

From the Baire category theorem, it follows that some $K_{j_0}$ has nonempty interior in $K$.
As $K^\circ$ is dense in $K$, this interior contains a nonempty
relatively open boundary patch $M_0\subset K^\circ$.
Thus the plaques at all points of $M_0$ have a common radius
$r_0=1/j_0$.

Now we need to choose a common transverse graph.
Noting that a convex boundary is locally a Lipschitz graph, we use standard
facts about convex functions and their supporting functionals
\cite[Sections 23--25]{Rock}.
By Rademacher's theorem, choose a differentiability point of such a graph in $M_0$. Translate this point to the origin and choose unitary coordinates $(w,z')$ with $w=x+it$ so that the tangent hyperplane to the graph is $\{x=0\}$ and its outward unit normal is $e_x$.

Then all outward unit supporting normals at nearby boundary points are
uniformly close to $e_x$.
Indeed, any subsequential limit of such normals supports $\Omega$ at $0$, where the outward unit normal is unique.
Consequently, in a sufficiently small chart contained in $M_0$, we have
\[
 \partial\Omega=\{x=-f(t,z')\},
 \qquad
 \Omega=\{x<-f(t,z')\}.
\]
Here $f$ is convex and Lipschitz on a real product neighborhood
and $f(0)=0$.

Write a supporting real functional at a point of a complex hyperplane
plaque as $\operatorname{Re}L$ with $L$ complex-linear.
At a relative interior point, $L$ annihilates the complex direction of
the plaque. That direction therefore equals $\ker L$ and is unique.
Put
\[
 g(t)=f(t,0),\quad a(t)=-g(t)+it,\quad q_t=(a(t),0).
\]
For a unit supporting normal at $q_t$, we can write the corresponding
complex-linear functional as
$L_t(w,z')=c_t w+\ell_t z'$.
The closeness of the normals to $e_x$ gives $|c_t|\ge1/2$ after
shrinking, while the rows $\ell_t$ remain bounded.
The plaque direction is $\ker L_t$, so its affine hyperplane has the
equation $w=a(t)+b(t)z'$, with $b(t)=-c_t^{-1}\ell_t$.
These rows are uniformly bounded, and $b(0)=0$ because the complex
tangent hyperplane at $0$ is $\{w=0\}$.

The uniform plaque radius now gives $R>0$ and an interval
$J$ about $0$ such that
\begin{equation}\label{eq:baire-plaques}
 P(t,z'):=\bigl(a(t)+b(t)z',z'\bigr)\in\partial\Omega,
 \quad t\in J,\quad |z'|<R,
\end{equation}
and all these points lie in the same boundary graph chart.
Here $J$ and $R$ are chosen small enough that the points
remain within distance $r_0$ of their centers and inside that chart.
We need to prove the following Claim.

\medskip
{\bf Claim:} The slope $b(t)$ is Lipschitz and constant.
\medskip

{\noindent\emph{Proof of Claim:}} For each fixed $|z'|<R$, the function
\[
 \phi_{z'}(t)=t+\operatorname{Im}\bigl(b(t)z'\bigr)
\]
is injective on $J$.
Indeed, if $\phi_{z'}(s)=\phi_{z'}(t)$, the two points in
\eqref{eq:baire-plaques} have the same graph coordinates
$(z',\operatorname{Im}w)$ and hence coincide.
They lie in the relative interiors of both plaques.
A supporting functional at this common point annihilates both complex
hyperplane directions, so the directions coincide. The two affine
hyperplanes, having a point in common, therefore coincide as well.
Their intersections with $\{z'=0\}$ give $a(s)=a(t)$, and taking
imaginary parts yields $s=t$.

Fix $s<t$ in $J$. The continuous function of $z'$
\[
 (t-s)+\operatorname{Im}\bigl((b(t)-b(s))z'\bigr)
\]
never vanishes on the connected ball $|z'|<R$ and is positive at $z'=0$.
Applying this also to $-z'$, we have
\[
 \left|\operatorname{Im}\bigl((b(t)-b(s))z'\bigr)\right|<t-s,
 \qquad |z'|<R.
\]
For a complex-linear row $\ell$,
$\sup_{|z'|<R}|\operatorname{Im}(\ell z')|=R|\ell|$.
Taking the supremum  yields that
\begin{equation}\label{eq:baire-slopes-Lipschitz}
 R|b(t)-b(s)|\le t-s,
\end{equation}
which implies that $b$ is Lipschitz on $J$.

Now choose $t\in J$ at which both $a$ and $b$ are differentiable.
This holds for almost every $t$, since both functions are Lipschitz.
A supporting real functional at $q_t$ annihilates its complex plaque
and hence has the form
\[
 \operatorname{Re}\bigl[c\bigl(w-a(t)-b(t)z'\bigr)\bigr]\le0,
 \quad (w,z')\in\overline\Omega,
 \quad c\in\mathbb C\setminus\{0\}.
\]
By substituting \eqref{eq:baire-plaques} for every $|z'|<R$ and $s\in J$,
we therefore obtain that
\[
 \operatorname{Re}\bigl[
 c\bigl(a(s)-a(t)+(b(s)-b(t))z'\bigr)\bigr]\le0
\]
with equality at $s=t$. Differentiating at this maximum gives
\[
 \operatorname{Re}\bigl[c\bigl(a'(t)+b'(t)z'\bigr)\bigr]=0,
 \quad |z'|<R.
\]
Subtracting the identity for $z'=0$ and varying $z'$ in a complex ball,
we obtain that $c\,b'(t)=0$ and therefore $b'(t)=0$.
Hence $b'=0$ almost everywhere. Since a Lipschitz function on an
interval is absolutely continuous, $b$ is constant on $J$.
As $b(0)=0$, we obtain $b\equiv0$.

Now by the Claim above, the plaques in \eqref{eq:baire-plaques} are given by
$w=a(t)=-g(t)+it$.
Uniqueness of the boundary graph therefore gives
$f(t,z')=g(t)$ for $t\in J$ and $|z'|<R$.
Shrink to a ball $B'$ and a rectangle $W$ inside the graph chart.
The boundary and the domain side are then respectively
$\Gamma\times B'$ and $G\times B'$, with $\Gamma$ and $G$ as stated.
Convexity and the Lipschitz property of $g$ follow from those of $f$.
\end{proof}

\begin{proof}[Proof of \Cref{thm:nonsmooth}]

Let $U$ be a ball on which $\Omega=D\cap U$ is convex. Work in a smaller ball
compactly contained in $U$, where $\partial\Omega=\partial D$.
By \Cref{lem:straight}, the rank hypothesis supplies an affine
complex hyperplane plaque through each point of the smaller true
boundary patch. Thus \Cref{prop:nonsmooth-product} applies.
Use the resulting coordinates, with
$W=(-h,h)+i(-\eta,\eta)$ and $g(0)=0$.

Choose $0<a<h$, and then $0<\eta_1<\eta$ so small that
$|g(t)|<a/2$ for $|t|\le\eta_1$. The planar domain
\[
 G_0=\{x+it:|t|<\eta_1,\ -a<x<-g(t)\}
\]
has closure contained in $W$. Its boundary consists of three line
segments and the Lipschitz graph arc $x=-g(t)$, so it is a Jordan
domain. A Riemann map $\phi:\D\to G_0$ extends homeomorphically to
the closures by Carath\'eodory's Jordan-domain theorem.
The preimage of the relative interior of the graph arc is a nonempty
open circle arc.

Fix $b\in\C^{n-1}\setminus\{0\}$ and choose $\epsilon>0$ such that
$sb\in B'$ for $|s|\le\epsilon$. Define
\[
 F(s,\lambda)=(\phi(\lambda),s b).
\]
This defines a holomorphic family with values in $D$, attached to the boundary along the same nonempty open arc. Since it varies with $s$, \Cref{prop:forbidden} yields a contradiction.
\end{proof}

\begin{proof}[Proof of \Cref{cor:disc-free}]
Suppose $\Omega=D\cap U$ is convex for a ball $U$ centered at a boundary
point. Every $q\in M=\partial D\cap U$ is locally convex. Choose a
smaller ball $V_q$ centered at $q$ with $\overline{V_q}\subset U$ and
use $D\cap V_q=\Omega\cap V_q$.

By \cite[Theorem 1.15]{PW}, it follows that no such $q$ can be a local holomorphic
peak point. Applying \Cref{lem:peak-criterion} to the bounded convex
domain $D\cap V_q$, we obtain a nonconstant analytic disc through $q$
in its boundary. Restricting the disc near its center so that its image
lies in $V_q$, it lies in the true boundary $M$.

Since the dimension is two, \Cref{lem:straight} yields an
affine complex disc at $q$ contained in $M$, even if the
disc was branched at its center. After rescaling, there is a map
$f_q:\D\to M$ with $f_q(0)=q$ and $\rank_\C\dd f_q(0)=1$.
Thus $M$ satisfies the hypotheses of \Cref{thm:nonsmooth}, which is a
contradiction.
\end{proof}

\section{Analytic tools}\label{sec:analytic-tools}

\subsection{\texorpdfstring{$C^1$}{C1} attached discs}

We use the following fixed-direction version of the $C^1$ Bishop construction in \cite[Section 3.1]{Suk}. The original construction is due to Bishop \cite{Bishop}.
\begin{lemma}\label{lem:bishop}
Let $E\subset\C^N$ be a $C^1$ maximally totally real submanifold through $0$, and let $\eta\in T_0E$ be fixed. Choose a smooth function $\psi\ge0$ on $\partial\D$, not identically zero, that vanishes on a nonempty open arc $I$. For every sufficiently small $\delta>0$, there is
\[
 u_\delta\in\OO(\D,\C^N)\cap C^0(\overline\D,\C^N)
\]
with boundary values
\begin{equation}\label{eq:bishop-boundary}
 u_\delta(\xi)=e_\delta(\xi)+i\delta\psi(\xi)\eta,
 \quad e_\delta(\xi)\in E,
\end{equation}
and $\|u_\delta\|_{C^0(\overline\D)}+\|e_\delta\|_{C^0(\partial\D)}=O(\delta)$.
\end{lemma}
\begin{proof}
After a complex-linear change of coordinates, write
$E=\{x+ih(x)\}$ locally, where $h\in C^1$, $h(0)=Dh(0)=0$,
and $\eta=a\in\R^N$.
A cutoff extension makes $h$ globally
Lipschitz with an arbitrarily small constant $L$; take $L<1/2$.
Indeed, on a ball of radius $2\varepsilon$,
$\sup|Dh|=o(1)$ and $\sup|h|=o(\varepsilon)$, while the cutoff
derivative is $O(\varepsilon^{-1})$.

Let $T$ be the circle Hilbert transform, with
$T(e^{ij\theta})=-i\operatorname{sgn}(j)e^{ij\theta}$ for $j\ne0$
and $T(1)=0$. Then $T^2=-\id+\operatorname{mean}$ and
$\|T\|_{L^2\to L^2}=1$. The Bishop equation
\[
 x=-T\bigl(h(x)+\delta\psi a\bigr)
\]
is a contraction on real-valued $L^2(\partial\D,\R^N)$ and has a
unique solution satisfying $\|x\|_2=O(\delta)$.
For the rotation difference
$\Delta_\theta x(\xi)=x(e^{i\theta}\xi)-x(\xi)$, commutation with $T$
gives
\[
 (1-L)\|\Delta_\theta x\|_2
 \le\delta|a|\,\|\Delta_\theta\psi\|_2
 \le C\delta|\theta|.
\]
The difference-quotient characterization of Sobolev spaces yields
$\|x\|_{W^{1,2}}=O(\delta)$. The same bound holds for
$y=h(x)+\delta\psi a$ by the Lipschitz composition rule.
Applying $T$ to the Bishop equation gives
$y=Tx+\operatorname{mean}(y)$, so $x+iy$ is the boundary value of a
holomorphic disc.

Sobolev embedding on the circle and the maximum
principle give its continuity on $\overline\D$ and the required
$O(\delta)$ bounds. For small $\delta$, the cutoff is inactive, so
$e_\delta=x+ih(x)$ lies in the original manifold. The inverse
complex-linear change of coordinates preserves the prescribed
displacement $i\delta\psi\eta$.
\end{proof}

\subsection{Extension along one disc}
We apply the local Baouendi--Tr\`eves approximation theorem \cite{BT} to finitely many $C^1$ CR matrix entries, all defined on the same larger neighborhood $U_S$.
Note the $C^2$ formulation with a fixed approximation neighborhood is recalled
in \cite[Section 2]{LNR}.

\begin{theorem}\label{thm:BT}
Let $S\subset\C^r$ be an embedded $C^2$ generic CR submanifold and
let $U_S$ be a fixed sufficiently small relatively open neighborhood of
$\zeta_0$ in $S$. There is a compact neighborhood
$K_{\mathrm{app}}\Subset U_S$ of $\zeta_0$, depending on $U_S$ but not
on the function, such that every $C^1$ CR function on $U_S$ is uniformly
approximable on $K_{\mathrm{app}}$ by holomorphic polynomials.
\end{theorem}

\begin{lemma}\label{lem:disc-extension}
Let $B:S\to\operatorname{Mat}_{a\times b}(\C)$ be $C^1$ and CR on
the fixed neighborhood $U_S$. Suppose
$\zeta\in\OO(\D,\C^r)\cap C^0(\overline\D,\C^r)$ satisfies
$\zeta(\partial\D)\subset K_{\mathrm{app}}$, where
$K_{\mathrm{app}}$ is as in \Cref{thm:BT}. Then there is
\[
 \mathcal B\in\OO(\D,\operatorname{Mat}_{a\times b}(\C))
 \cap C^0(\overline\D,\operatorname{Mat}_{a\times b}(\C))
\]
with $\mathcal B|_{\partial\D}=B\circ\zeta$. If $B\circ\zeta$ is sufficiently close to a fixed matrix $B_0$ on the circle, then $\mathcal B$ is close to $B_0$ on the entire closed disc.
In particular, a fixed invertible square block of $B_0$ remains invertible there for sufficiently small oscillation.
\end{lemma}
\begin{proof}
Approximate the finitely many entries simultaneously by polynomial matrices $P_\nu$. The maximum principle applied entrywise gives
\[
 \sup_{\overline\D}|P_\nu\circ\zeta-P_\mu\circ\zeta|
 \le C\sup_{K_{\mathrm{app}}}|P_\nu-P_\mu|\longrightarrow0.
\]
Their uniform limit is the required matrix $\mathcal B$.
The last assertion follows by applying the maximum principle to $\mathcal B-B_0$ and then using the continuity of matrix inversion.
\end{proof}

\subsection{Real convex geometry}\label{sec:real-geometry}

Let $M$ be a locally convex $C^2$ hypersurface. Its outward unit
normal $\nu$ is $C^1$. We use the sign convention
\[
 A=\dd\nu:TM\to TM,\quad A\ge0.
\]
Put
\[
 H=TM\cap JTM,\quad N=\ker A.
\]
Thus $H$ is the real $J$-invariant complex tangent bundle. For an
interior-negative $C^2$ defining function, the Levi quadratic form on $H$
is a positive scalar multiple of
\begin{equation}\label{eq:real-Levi}
 \langle Av,v\rangle+\langle AJv,Jv\rangle.
\end{equation}
Consequently convexity implies
\begin{equation}\label{eq:K}
 K=\ker\Levi=N\cap JN
   =\{v\in H:Av=AJv=0\}.
\end{equation}
Here $K$ denotes the real realization of the Levi kernel, identified
with a complex subspace of $\C^n$. The identity uses positivity:
$\langle Av,v\rangle=0$ implies $Av=0$ for a positive semidefinite
self-adjoint operator. In general $N\cap H$ need not be $J$-invariant,
and it need not equal $K$.

\begin{lemma}\label{lem:no-strongly-pseudoconvex}
A locally convex $C^2$ boundary patch of $D$ contains no strongly
pseudoconvex point.
\end{lemma}
\begin{proof}
If $p$ were such a point, no affine complex disc through $p$ could be locally contained in the boundary. Indeed, for a nonzero direction $v$ of such a disc and a $C^2$ defining function $\rho$, differentiating $\rho(p+\lambda v)=0$ first shows that $v$ is complex tangent and then gives
\[
 \sum_{j,\ell}\rho_{z_j\bar z_\ell}(p)
 v_j\overline{v_\ell}=0,
\]
contrary to strong pseudoconvexity. Applying \Cref{lem:peak-criterion} to a sufficiently small bounded convex local intersection, we would obtain a local holomorphic peak function at $p$, contradicting \cite[Theorem 1.15]{PW}.
\end{proof}

For the positive-rank construction below, assume the Levi rank
is constantly $r$ with $1\le r<n-1$, and put $k=n-1-r$.
After choosing a point of maximal real rank and shrinking, we may
also assume that $m=\rank_\R A$ is constant.
The Levi-flat case will be excluded separately using
\Cref{thm:nonsmooth} in the proof of \Cref{cl:main}.

\begin{lemma}\label{lem:Legendre}
Let $f\in C^2(V)$ be convex on an open convex set $V\subset\R^d$,
with $\rank D^2f\equiv m$. Near each point, its gradient image is an
embedded $C^2$ submanifold of dimension $m$. Its local constant-gradient
fibers are relatively open subsets of affine $(d-m)$-planes.
The local constant-normal fibers of the graph of $f$ are also affine
pieces, with tangent space equal to the kernel of its shape operator.
\end{lemma}
\begin{proof}
The case $m=0$ is immediate. Otherwise choose real orthogonal coordinates $(x,y)\in\R^m\times\R^{d-m}$ with $f_{xx}>0$ on a small neighborhood. Set $p=f_x(x,y)$, and solve $x=X(p,y)$ by the $C^1$ implicit function theorem. Define
\[
 g(p,y)=p\cdot X(p,y)-f(X(p,y),y).
\]
The cancellation identities
\begin{equation}\label{eq:Legendre-first}
 g_p=X,\qquad g_y=-f_y(X,y)
\end{equation}
show that $g\in C^2$, although $X$ is only known to be $C^1$. Constant Hessian rank gives the vanishing Schur complement
\[
 f_{yy}-f_{yx}f_{xx}^{-1}f_{xy}=0.
\]
Differentiating $f_x(X(p,y),y)=p$ gives
$X_y=-f_{xx}^{-1}f_{xy}$. Hence
\[
 g_{yy}=-f_{yy}+f_{yx}f_{xx}^{-1}f_{xy}=0.
\]
Choose a product neighborhood in the $(p,y)$ variables, translating
the basepoint in the $y$ variables to $0$. On this neighborhood,
\[
 g(p,y)=b(p)+a(p)\cdot y.
\]
Both $a$ and $b$ are $C^2$: recover $b$ by evaluating at $y=0$ and each component of $a$ by taking a fixed difference between $y=0$ and $y=he_j$. Hence
\begin{equation}\label{eq:gradient-image}
 \nabla f(X(p,y),y)=(p,-a(p)),
\end{equation}
an embedded $C^2$ graph in $p$.

For fixed $p$, $X=g_p=Db(p)+(Da(p))^{\mathsf T}y$ is affine in $y$, and
\[
 f(X(p,y),y)=p\cdot X(p,y)-b(p)-a(p)\cdot y
\]
is affine in $y$ as well.
The $y$ coordinates remain free, so these parametrizations have
rank $d-m$ and their images are relatively open in affine
$(d-m)$-planes.
In \eqref{eq:gradient-image}, $p$ is part of the gradient itself.
Therefore these are exactly the local constant-gradient fibers.
The normal of a convex graph is a smooth injective function of its
gradient, so its fibers are the same.

Their tangent spaces coincide with the kernel of the differential of the Gauss map, and hence with the kernel of the shape operator. The second fundamental form is, up to a sign depending on the orientation, a nonzero scalar multiple of $D^2f$. Since the induced metric is invertible, the shape operator has rank $m$.
\end{proof}

\begin{remark}\label{rem:image-map}
The preceding lemma also gives a $C^2$ structure on the local image
$\Sigma$ of the Gauss map of a constant-real-rank convex $C^2$
hypersurface. To see this, after applying a real rigid motion, represent the hypersurface as the graph of a convex function. The map
\[
 u\longmapsto\frac{(-u,1)}{\sqrt{1+|u|^2}}
\]
is a smooth diffeomorphism onto a hemisphere. It carries the $C^2$
gradient image in \Cref{lem:Legendre} onto $\Sigma$, up to the choice
of orientation. Moreover,
$T_{\nu(p)}\Sigma=\im A_p$, and the original map $M\to\Sigma$
is a $C^1$ submersion. Thus the image has $C^2$ regularity even though
the Gauss map is only $C^1$.
\end{remark}

\subsection{The direct affine-fiber construction}\label{sec:easy}

Suppose at some point $p$ there is $v\in N_p\setminus H_p$. Choose $0\ne b\in K_p$. Since $Jv$ is not tangent to $M$, change the sign of $v$ to have $\dd\rho_p(Jv)<0$ for an interior-negative local defining function $\rho$.
The local constant-normal fiber through $p$ is an open affine piece with tangent $N_p$, by \Cref{lem:Legendre}. Since $b,Jb,v\in N_p$, the points
\[
 p+xv+sb\quad (x\in\R,\ s\in\C)
\]
lie in $M$ for all sufficiently small $x,s$, where
$sb=(\Rea s)b+(\Ima s)Jb$.

Choose a smaller closed product of parameter intervals inside this
affine piece and the defining-function chart.
Since $\dd\rho_p(Jv)<0$, continuity and a further shrinking give
constants $c,y_0>0$ such that
\[
 \rho(p+xv+sb+yJv)\le-cy
 \quad(0\le y\le y_0)
\]
uniformly for those parameters.
This follows by integrating the uniformly negative derivative in
the fixed $Jv$ direction from the boundary point $p+xv+sb$.
Thus the map
\begin{equation}\label{eq:easy-family}
 G(s,w)=p+wv+sb
\end{equation}
is holomorphic into the local convex side for $s$ small and $w$ in a
small upper rectangle. A Riemann map from $\D$ to that rectangle extends homeomorphically to the closures because the rectangle is a Jordan domain.

The preimage of the relative interior of the lower side is a
nonempty open circle arc. Composing \eqref{eq:easy-family} with this map gives a nontrivial
holomorphic family attached to the true boundary on that arc. Its
slices are distinct because $b\ne0$.
Since $N$ is a continuous bundle on a constant-rank patch, the condition $N\not\subset H$ is open. If it occurs, the direct construction suffices.
It remains to treat a patch on which
\begin{equation}\label{eq:difficult-branch}
 N\subset H\quad\text{at every point}.
\end{equation}

\section{The conormal reduction for the remaining case}\label{sec:conormal}

We now work on a constant-Levi-rank, constant-real-rank patch satisfying
$N\subset H$, with $1\le r<n-1$ and $k=n-1-r$.
Our goal is to encode each affine complex fiber in terms of its direction and position, to obtain a $C^1$ maximally totally real parameter space for the fibers, and to retain a linear formula for recovering points in the original $\C^n$.

The regularity needed for the CR coefficients comes from a separate
$C^2$ conormal parameter manifold.
We can translate the basepoint of $M$ to $0$ in the original $\C^n$ and center the conormal parameter coordinates at $\zeta_0=0$.
The quotient basepoint will then be $(0,B(0)0)=(0,0)$, as required
for the Bishop construction.

\subsection{The projective image}
For a real vector $\nu\in\C^n$, define the complex-linear covector
\[
 \theta_\nu(v)=\sum_j\overline{\nu_j}v_j
 =\langle\nu,v\rangle+i\langle J\nu,v\rangle.
\]
The projective conormal map is
$a(p)=[\theta_{\nu(p)}]$, valued in the projective space of complex
covectors.
Let $\nu(p)$ be the outward unit normal vector of $M$ at $p$. We may assume that $c(p)=\theta_{\nu(p)}(e_n)\ne0$, and represent $a(p)$ by
\[
 \alpha(p)=\theta_{\nu(p)}/c(p),\qquad\alpha_n=1.
\]
Equivalently, $\alpha_j(p)=\overline{\nu_j(p)}/\overline{\nu_n(p)}$.
 Note that we can identify the affine chart
$\{\alpha_n=1\}$ with $\C^{n-1}$ by omitting the last component.

\begin{lemma}\label{lem:annihilator}
At each point,
\begin{equation}\label{eq:annihilator}
 \Ann_\C\Span_\C\{\alpha,\dd\alpha(TM)\}=K.
\end{equation}
In particular,
\[
 \dim_\C\Span_\C\dd\alpha(TM)=r.
\]
\end{lemma}

\begin{proof}
If $\alpha(v)=0$, then $v,Jv\in TM$. For each $X\in TM$, we have
\[
 \dd\alpha(X)(v)=c^{-1}\theta_{AX}(v)-c^{-2}\theta_{\nu}(v)dc(X)=c^{-1}\theta_{AX}(v).
\]
The last equality follows from the fact that $\theta_\nu(v)=0$. If $d\alpha(X)(v)=0$ for each $X\in TM$, then
\[
\langle AX,v\rangle=0,\ \langle JAX,v\rangle=0
\]
from the definition of $\theta_{\nu}$. Since $\langle JAX,v\rangle=-\langle AX,Jv\rangle$, both $v$ and $Jv$ are orthogonal to $\im A$. Self-adjointness of $A$ then gives $Av=AJv=0$, so $v\in K$ by \eqref{eq:K}. The converse follows from the same calculation.

Note that the span including $\alpha$ has complex dimension $n-k=r+1$. Since $\alpha_n=1$ while $(\dd\alpha)_n=0$, the complex span of the differential alone has dimension exactly $r$.
\end{proof}

Now assume that \eqref{eq:difficult-branch} holds. Then $\langle J\nu,u\rangle=0$ for every $u\in N$, and hence $J\nu\perp N$. Thus
\begin{equation}\label{eq:Jnu-image}
 J\nu\in N^\perp\cap TM=\im A.
\end{equation}
On the unit sphere, the kernel of the differential of
$\nu\mapsto[\theta_\nu]$ is exactly $\R J\nu$.
Indeed, a tangent vector $w$ is in this kernel precisely when
$\theta_w$ is a complex multiple of $\theta_\nu$, so
$w\in\Span_\R\{\nu,J\nu\}$.
As $w\perp \nu$, we have $w\in\mathbb{R}J\nu$.
Since $T_\nu\Sigma=\im A$ contains $J\nu$,
projectivization restricted to $\Sigma$ has real rank $m-1$.
The Gauss map $M\to\Sigma$ is a submersion, so
\[
 \rank_\R\dd a=m-1.
\]
By \Cref{rem:image-map}, the unit normal image $\Sigma$ is $C^2$.
Since the projectivization restricted to $\Sigma$ is a $C^2$ map of constant real rank $m-1$, the constant-rank theorem yields a locally embedded $C^2$ image $S'$ in the affine conormal chart.
Its real dimension is $m-1$. By \Cref{lem:annihilator}, the complex span of its real tangent space has dimension $r$.

For any real subspace $V$ of a complex vector
space,
\[
 \dim_\C(V\cap iV)=\dim_\R V-\dim_\C\Span_\C V.
\]
Consequently $S'$ has constant CR dimension and CR codimension
\begin{equation}\label{eq:q-c}
 d_{\mathrm{CR}}=m-1-r,
 \quad c_{\mathrm{CR}}=2r+1-m,
\end{equation}
respectively. Their nonnegativity gives $r+1\le m\le2r+1$.

Now we write $W_\alpha=\Span_\C T_\alpha S'$ and choose a fixed complex-linear
projection $\pi:\C^{n-1}\to\C^r$ that is injective on $W_{\alpha_0}$.
For clarity, the nearby spaces $W_\alpha$ need not be parallel to
$W_{\alpha_0}$. They vary continuously and have constant dimension
$r$: choose a continuous local complex frame $w_1(\alpha),\ldots,w_r(\alpha)$.
The determinant of the columns $\pi w_1(\alpha),\ldots,\pi w_r(\alpha)$ is
nonzero at the basepoint, hence remains nonzero nearby. Thus the
same fixed projection is injective on every $W_\alpha$ after shrinking.
In particular, $\pi|_{S'}$ is an immersion and, on a smaller
neighborhood, a $C^2$ embedding. We may center its image by setting
\[
 S=\pi(S')-\pi(\alpha_0)\subset\C^r.
\]
Since $\Span_\C TS=\C^r$, $S$ is a $C^2$ generic CR
manifold. Denote its inverse parametrization, with the last component
$1$ restored, by $\alpha:S\to\C^n$.

Let $u,iu\in TS$ and $U,V\in TS'$
be their respective lifts. Then $\pi(V-iU)=0$. Injectivity on the
entire complex span of $TS'$ gives $V=iU$, which is precisely
the CR condition for the inverse. Thus the parameter map
$\zeta:M\to S$ is obtained by taking
$\pi(\alpha(p))-\pi(\alpha_0)$ for the normalized conormal at $p$,
so $\zeta(0)=0$. It is a $C^1$
submersion, and the normalized conormal at $p$ is
$\alpha(\zeta(p))$. Note that the map $\alpha:S\to\C^n$ is $C^2$ and CR.

\subsection{The full matrix of first jets}
The first derivatives of the normalized conormal have the following
CR decomposition.
\begin{lemma}\label{lem:jets}
There are unique $C^1$ row-vector functions $A_1,\ldots,A_r$ on $S$ satisfying
\begin{equation}\label{eq:jet-decomp}
 \dd\alpha=\sum_{j=1}^r A_j\,\dd\zeta_j\quad\text{on }TS.
\end{equation}
Every entry of every $A_j$ is CR. The matrix
\begin{equation}\label{eq:B}
 B(\zeta)=\begin{pmatrix}\alpha(\zeta)\\A_1(\zeta)\\\vdots\\A_r(\zeta)\end{pmatrix}
\end{equation}
has rank $r+1$, and $\ker B(\zeta(p))=K_p$.
\end{lemma}
\begin{proof}
Because $S$ is generic in $\C^r$, the restrictions of the $r$ forms
$\dd\zeta_j$ are linearly independent in the complexified cotangent
bundle. As the common kernel of $\dd\zeta_j$ is $T^{0,1}S$, they form a basis of its
annihilator. Since $\alpha$ is $C^2$ and CR, this gives
\eqref{eq:jet-decomp} with unique coefficients $A_j$.

In any $C^2$ coordinate chart on $S$, these coefficients are obtained by solving a
linear system whose coefficients and right-hand side are $C^1$ and
whose chosen minor stays invertible. Thus $A_j\in C^1$.
Taking exterior derivatives, we have
\[
 0=\sum_j\dd A_j\wedge\dd\zeta_j.
\]
Now evaluate on $L\in T^{0,1}S$ and an arbitrary complexified tangent vector $X$. Since $\dd\zeta_j(L)=0$, we obtain
\[
 \sum_j(LA_j)\,\dd\zeta_j(X)=0.
\]
Note that genericity means that the map
\[
X\mapsto(\dd\zeta_1(X),\ldots,\dd\zeta_r(X))
\]
is surjective onto $\C^r$. This implies that $LA_j=0$ for every $j$.
Finally, the row span of \eqref{eq:B} is the span in \Cref{lem:annihilator}, since $\zeta$ is a submersion. This completes the proofs of both the rank and kernel assertions.
\end{proof}

\subsection{The totally real quotient}
In this subsection, we construct a totally real space on which Lemma \ref{lem:bishop} can be applied. We first retain both the conormal parameter and the original point by
forming the graph
\[
 \iota:M\longrightarrow S\times\C^n,\quad
 \iota(p)=(\zeta(p),p),\quad \Lambda=\iota(M).
\]
The projection to the second factor identifies $\Lambda$ with $M$.
The purpose of this lift is to reduce its complex tangent directions
to the Levi kernel, as the next lemma shows.

\begin{lemma}\label{lem:lift}
The $C^1$ lift satisfies
\[
 T\Lambda\cap iT\Lambda=\{0\}\times K.
\]
\end{lemma}
\begin{proof}
Suppose the lifted tangent of $X\in TM$ and its multiplication by $i$ both belong to $T\Lambda$. Then $JX\in TM$ and
\[
 \dd\zeta(JX)=i\dd\zeta(X).
\]
Since the inverse parametrization $\alpha$ is CR, the same relation holds for the complete normalized conormal.
For a $C^2$ defining function $\rho$, we have $\nu=\frac{\nabla\rho}{|\nabla\rho|}$ and $\theta_{\nu}=2\frac{\partial\rho}{|\nabla\rho|}$. This implies  $\alpha_j=\rho_{z_j}/\rho_{z_n}$ with $j<n$.
Then
\[
\dd\rho(X)=\sum_l\rho_{z_l}X_l+\sum_l\rho_{\bar{z}_l}\overline{X}_l,\ \dd\rho(JX)=i\sum_l\rho_{z_l}X_l-i\sum_l\rho_{\overline{z}_l}\overline{X}_l.
\]
After a direct differentiation, we get
\begin{align*}
 \dd\alpha_j(JX)-i\dd\alpha_j(X)
 =-\frac{2i}{\rho_{z_n}}
 \sum_{\ell=1}^n
 \bigl(\rho_{z_j\bar z_\ell}-\alpha_j\rho_{z_n\bar z_\ell}\bigr)
 \overline{X_\ell}.
\end{align*}
By setting $h_j=\sum_\ell\rho_{z_j\bar z_\ell}\overline{X_\ell}$,
the identity $\dd\alpha_j(JX)=i\dd\alpha_j(X)$ gives that
$$
h_j=\alpha_jh_n=(\rho_{z_j}/\rho_{z_n})h_n
$$
for $j<n$,
and the same equality is automatic for $j=n$.
Since $X$ is complex tangent, we have
\[
 \Levi_\rho(X)=\sum_j X_jh_j
 =\frac{h_n}{\rho_{z_n}}\sum_j\rho_{z_j}X_j=0.
\]
By positivity and \eqref{eq:real-Levi}, $AX=AJX=0$.
Thus $X\in K$, and the normal derivative vanishes, giving
$\dd\zeta(X)=0$.

Conversely, if $X\in K$, the normal derivatives in both the $X$
and $JX$ directions vanish. Hence both corresponding
$\zeta$ derivatives vanish, and the lifted vectors are
$(0,X)$ and $(0,JX)$.
\end{proof}

Now we consider the map
\[
 Q:S\times\C^n\to\C^r\times\C^{r+1},
 \quad Q(\zeta,z)=(\zeta,B(\zeta)z).
\]
We use standard facts about CR functions on generic submanifolds, and one can see \cite{BER} for more details. Note that every $C^1$ CR function on a generic submanifold has a unique complex-linear extension of its real differential to the ambient complex vector space.

Indeed, write an ambient vector as $X+iY$ with $X,Y\in TS$ and define the extension to be $df(X)+i\,df(Y)$. Independence of this decomposition is exactly the CR identity on $TS\cap iTS$ and the uniqueness follows from $TS+iTS=\C^r$.
We denote the pointwise complex-linear extension of the differential by $\mathscr DB$.
Then the differential of $Q$ on $S\times\mathbb{C}^n$ is the restriction of the complex-linear map
\begin{equation}\label{eq:formalDQ}
 \mathscr DQ(u,w)=\bigl(u,(\mathscr DB(u))z+B(\zeta)w\bigr),
 \quad\ker_\C\mathscr DQ=\{0\}\times K.
\end{equation}

\begin{lemma}\label{lem:quotient}
Locally $q=Q|_\Lambda$ has a $C^1$ embedded image $E$ that is maximally totally real in $\C^{2r+1}$. Its real differential has kernel $\{0\}\times K$ on $T\Lambda$.
\end{lemma}
\begin{proof}
At a point $(\zeta(p),p)$ of $\Lambda$, the complex kernel of
\eqref{eq:formalDQ} is $\{0\}\times K_p$, and \Cref{lem:lift} shows
that this entire kernel lies in $T\Lambda$. Hence the real kernel of
$\dd q$ is exactly $\{0\}\times K_p$, and $q$ has constant real rank
\[
 (2n-1)-2k=2r+1.
\]
The $C^1$ constant-rank theorem supplies that $E$ is an embedded submanifold with $\operatorname{dim}_{\mathbb{R}}E=2r+1$ and makes $q$ a submersion onto $E$.

To check total reality, fix $p\in M$ and suppose
$e,ie\in T_{q(p)}E$. As $q$ is a submersion, we can choose $X,Y$ in the
\emph{same} tangent space $T_{(\zeta(p),p)}\Lambda$ with
$\dd Q(X)=e$ and $\dd Q(Y)=ie$. By the complex linearity of
\eqref{eq:formalDQ},
$Y-iX\in\{0\}\times K_p\subset T_{(\zeta(p),p)}\Lambda$.
Consequently, $iX$ belongs to this tangent space.
By \Cref{lem:lift}, we have $X\in\{0\}\times K_p$ and hence $e=0$. Thus $TE\cap iTE=\{0\}$, and the dimension calculation shows that $E$ is maximally totally real.
\end{proof}

We identify $M$ with $\Lambda$ via $\iota$ and use the same letter
$q$ for the resulting quotient map on $M$. Explicitly, define
\begin{equation}\label{eq:quotient-coordinates}
 \beta(p)=B(\zeta(p))p,\quad
 q(p)=(\zeta(p),\beta(p)),\quad E=q(M).
\end{equation}
The role of the two coordinate groups is expressed by
\begin{equation}\label{eq:affine-fiber-equation}
 B(\zeta(p))z=\beta(p)
 \quad\Longleftrightarrow\quad z\in p+K_p.
\end{equation}
Thus $\zeta$ determines the direction matrix, and $\beta$ fixes the
position of the affine fiber. The quotient removes its $k$ complex
degrees of freedom. The next lemma verifies that these affine planes
are actual boundary fibers with a uniform local size.

\subsection{Uniform affine fibers}

Choose a fixed set $\mathcal I$ of $r+1$ column indices such that
$B_{\mathcal I}$ is invertible near the basepoint.
Let $\mathcal J$ be the complementary set of $k$ column indices.

\begin{lemma}\label{lem:fibers}
After translating the basepoint in $\C^n$ to $0$, there are a fixed neighborhood $E_0\subset E$ and a fixed ball $\mathbb{B}_\epsilon^k$ such that
$E_0\times \mathbb{B}_\epsilon^k$
is a $C^1$ local parametrization of $M$.
\end{lemma}
\begin{proof}
We now consider the map
\[
\Theta:M\to E\times \mathbb{C}^k,\ \Theta(p)=(q(p),p_{\mathcal{J}})
\]
with
\[
q(p)=(\zeta(p),B(\zeta(p))p).
\]
If $d\Theta_p(X)=0$, then $X\in\operatorname{ker}B(\zeta(p))=K_p$. Since $X_{\mathcal{J}}=0$, we have $B_{\mathcal{I}}X_{\mathcal{I}}=0$, which means that $X_{\mathcal{I}}=0$ and hence $X=0$.

Since $\operatorname{dim}_{\mathbb{R}}M=\operatorname{dim}_{\mathbb{R}}(E\times\mathbb{C}^k)=2n-1$, the injective linear map $d\Theta_p$ is an isomorphism.  By the \(C^1\) inverse function theorem, after shrinking around
\(p = 0\), it follows that there are a relatively open neighborhood \(E_0\) of \(q(0)\) in \(E\), a ball \(\mathbb{B}_\epsilon^k\), and an open set
\(U \subset M\) such that
\[
\Theta|_U : U \longrightarrow E_0 \times B_\epsilon^k
\]
is a \(C^1\) diffeomorphism. Denote its inverse by
$R : E_0 \times B_\epsilon^k \longrightarrow U$.
Then \(\Theta \circ R = \mathrm{id}_{E_0 \times B_\epsilon^k}\) and \(R \circ \Theta|_U = \mathrm{id}_U\).

For $e=(\zeta,\beta)\in E_0$ and $s\in \mathbb{B}_{\epsilon}^k$,
the definition of $\Theta$ gives
\begin{equation}\label{eq:R}
 R_{\mathcal J}(e,s)=s,\quad
 R_{\mathcal I}(e,s)=B_{\mathcal I}(\zeta)^{-1}\bigl(\beta-B_{\mathcal J}(\zeta)s\bigr).
\end{equation}
\end{proof}

Set \(e_0 = q(0) = (0,0)\), \(\zeta_0 = \zeta(0)\), and \(\nu_0 = \nu(0)\). All shape operators and differentials in this paragraph are evaluated at 0.

By \eqref{eq:Jnu-image}, we can solve $Av=J\nu_0$.
Since $Av\ne0$ and $A$ is positive
semidefinite, we have $\langle Av,v\rangle>0$.
Noting that
$\nabla\rho=|\nabla\rho|\nu$ on $M$, we have
\[
 \dd\rho(Jv)
 =-|\nabla\rho|\langle J\nu_0,v\rangle
 =-|\nabla\rho|\langle Av,v\rangle<0.
\]
Moreover, $\theta_{Av}=\theta_{J\nu}=-i\theta_\nu$.
Differentiating $\alpha=\theta_\nu/c$ along $v$, the numerator
and denominator both acquire the factor $-i$, and therefore
$\dd\alpha(v)=0$.
Since $\zeta=\pi\circ\alpha$ up to a translation, we get
$\dd\zeta(v)=0$. Differentiating $q(p)=(\zeta(p),B(\zeta(p))p)$ at $p=0$ gives
\begin{equation}\label{eq:vertical}
 \eta_0=\dd q(v)=(0,b_0),\quad b_0=B(\zeta_0)v.
\end{equation}
Note that the vector $b_0$ is nonzero. If not, then $v\in\ker B=K$ would imply
$Av=0$, which contradicts $Av=J\nu_0$.

For $e=(\zeta,\beta)$, we can define
\[
 W_{\mathcal J}(e)=0,\quad W_{\mathcal I}(e)=B_{\mathcal I}(\zeta)^{-1}b_0.
\]
Then \(B(\zeta)W(e) = b_0\) for every \(e = (\zeta, \beta) \in E_0\). In particular, \(B(\zeta_0)(W(e_0) - v) = 0\), so \(W(e_0) - v \in K_0\).
Since \(K_0\) is complex and contained in \(T_0M\), we have \(i(W(e_0) - v) \in T_0M\), which implies that
\[
d\rho_0(iW(e_0)) = d\rho_0(Jv) < 0.
\]

For small $s,t$ and $e\in E_0$, we denote
 \[
 \Phi(e,s,t)=\dd\rho_{R(e,s)+itW(e)}(iW(e)),
 \]
which satisfies $\Phi(e_0,0,0)<0$. After shrinking $E_0$ and $B_{\epsilon}^k$, we can choose a smaller product whose closure is compactly contained in the product supplied by \Cref{lem:fibers}. Continuity of $d\rho$ and of $W$ then makes the derivative in the $iW$ direction uniformly negative on that closed product and on a short surrounding tube. Since $R(e,s)\in M$, we obtain
\begin{equation}\label{eq:uniform-inward}
 \rho\bigl(R(e,s)+itW(e)\bigr)=\int_0^t\Phi(e,s,\tau)d\tau\leq-ct
\end{equation}
for each $e\in E_0,\ |s|<\epsilon,\ 0\le t\le t_0$. Here, the constant $c$ satisfies $\Phi\leq -c$.

\section{Proof of \Cref{cl:main}}\label{sec:family}

\begin{proposition}\label{prop:family}
Let $\Omega\subset\C^n$ be a convex domain and
$M\subset\partial\Omega$ be a nonempty relatively open $C^2$
hypersurface patch of constant Levi rank $r$, where $1\le r<n-1$.
There exist $\epsilon>0$, a nonempty open arc $I\subset\partial\D$,
and a holomorphic map
\[
 F:\D_\epsilon\times\D\longrightarrow\Omega
\]
that extends continuously to $\D_\epsilon\times\overline\D$.
Its slices are distinct, and $F(s,I)\subset M$ for all $|s|<\epsilon$.
\end{proposition}
\begin{proof}
Choose the largest real rank $m$ attained by the shape operator
on the constant-Levi-rank patch.
It is positive, since positive Levi rank rules out $A=0$.
A nonzero $m\times m$ minor persists near a point attaining this rank,
and maximality makes the real rank constant on a smaller open patch.
If $N\not\subset H$ at any point of the resulting patch, the direct
construction in \Cref{sec:easy} gives the required family.
Otherwise $N\subset H$ throughout the patch, and we use
$S,B,E,R,\eta_0$ from \Cref{sec:conormal}.

We first construct a Bishop disc in the quotient coordinates.
Fix the product neighborhood and the constants in
\eqref{eq:uniform-inward}. Choose a sufficiently small relatively
open neighborhood $U_S$ on which all entries of $B$ are $C^1$ and CR,
and then choose $K_{\mathrm{app}}$ from \Cref{thm:BT}. These choices
precede the choice of the Bishop-disc scale.

Apply \Cref{lem:bishop} to $E$ at $0$, with the fixed tangent vector
$\eta_0=(0,b_0)$ from \eqref{eq:vertical}. Choose a smooth
$\psi\ge0$, not identically zero, that vanishes on a nonempty open
arc $I\subset\partial\D$. For small $\delta>0$, the resulting disc
satisfies
\[
 u_\delta=(\zeta_\delta,\beta_\delta),\quad
 u_\delta(\xi)=e_\delta(\xi)+i\delta\psi(\xi)(0,b_0),
 \quad e_\delta(\xi)\in E.
\]
Writing $e_\delta=(\widehat\zeta_\delta,
\widehat\beta_\delta)$ on the circle makes the two components explicit.
\[
 \zeta_\delta=\widehat\zeta_\delta,\quad
 \beta_\delta=\widehat\beta_\delta+i\delta\psi b_0.
\]
Since $E\subset S\times\C^{r+1}$, the first identity gives
\begin{equation}\label{eq:whole-circle}
 \zeta_\delta(\partial\D)\subset S.
\end{equation}
The $O(\delta)$ estimate and the fact that $K_{\mathrm{app}}$
contains a relative neighborhood of $0$ in $S$ ensure that this
circle image lies in $K_{\mathrm{app}}$ for small $\delta$.

\Cref{lem:disc-extension} implies that there is a matrix
$\mathcal B_\delta$, holomorphic on $\D$ and continuous on
$\overline\D$, whose boundary values are
\[
 \mathcal B_\delta(\xi)=B(\zeta_\delta(\xi)),
 \quad \xi\in\partial\D.
\]
As $\delta\to0$, these boundary values
converge uniformly to $B(0)$, and the maximum principle gives the
same convergence on $\overline\D$. In particular,
\[
 \sup_{\overline\D}
 \left\|B_{\mathcal I}(0)^{-1}
 \bigl(\mathcal B_{\delta,\mathcal I}-B_{\mathcal I}(0)\bigr)
 \right\|<\tfrac12
\]
for small $\delta$. Hence $\mathcal B_{\delta,\mathcal I}$ is
invertible on the entire closed disc.

Now we seek $z$ satisfying the holomorphic linear system
\[
 \mathcal B_\delta(\lambda)z=\beta_\delta(\lambda).
\]
Its boundary values are the affine-fiber equations
\eqref{eq:affine-fiber-equation}, with the prescribed vertical
displacement in $\beta$. Choosing the free coordinates $z_{\mathcal J}=s$
determines the remaining coordinates uniquely.
We first construct the family for $s$ in a small ball in $\mathbb{C}^k$.
Define, jointly in
$s\in\C^k$ and $\lambda\in\D$,
\begin{equation}\label{eq:F}
 F_{\mathcal J}(s,\lambda)=s,\quad
 F_{\mathcal I}(s,\lambda)=\mathcal B_{\delta,\mathcal I}(\lambda)^{-1}
 \bigl(\beta_\delta(\lambda)-\mathcal B_{\delta,\mathcal J}(\lambda)s\bigr).
\end{equation}
The inverse matrix is holomorphic, so $F$ is jointly holomorphic and
continuous up to the circle. On the circle, \eqref{eq:R} and the
linear dependence on $\beta$ give the exact identity
\[
 F(s,\xi)=R(e_\delta(\xi),s)+i\delta\psi(\xi)W(e_\delta(\xi)).
\]

Fix a smaller ball of $s$-parameters with closure in
$\mathbb{B}_\epsilon^k$. Then choose $\delta$ so small that
$e_\delta(\partial\D)\subset E_0$ and
$\delta\|\psi\|_\infty<t_0$. Then \eqref{eq:uniform-inward} yields
\begin{equation}\label{eq:boundary-F}
 \rho(F(s,\xi))\le-c\delta\psi(\xi).
\end{equation}
The reconstruction formula gives a uniform bound
$|F(s,\xi)|\le C_0(\delta+|s|)$ on the circle, because the basepoint
is $0$, $\beta_\delta=O(\delta)$, and the matrix inverse is uniformly
bounded.
Choose a Euclidean ball $V=B(0,a)$ whose closure is contained in
the defining-function chart.

Shrink the parameter ball, independently of $\delta$, so that
$C_0|s|<a/4$. Then decrease $\delta$ so that
$C_0\delta<a/4$ and all earlier requirements hold.
Every boundary value then lies in $B(0,a/2)$.
The intersection $\Omega_0=\Omega\cap V$ is convex.
By \eqref{eq:boundary-F},
$F(s,\partial\D)\subset\overline\Omega_0$, and the nonempty open
set of circle points where $\psi>0$ has image in $\Omega_0$.

For each fixed $s$, the real coordinate functions of
$F(s,\cdot)$ are harmonic and continuous on $\overline\D$. Each
interior value is therefore a Poisson average of boundary values and
belongs to the closed convex set $\overline\Omega_0$. If an interior
value were on $\partial\Omega_0$, a nonzero real supporting affine functional
there would attain its maximum on the disc. The harmonic maximum
principle would make it constant. Continuity to the circle would then
contradict the boundary values in $\Omega_0$, where every nonzero supporting
functional is strictly smaller.
Thus
\[
 F(s,\D)\subset\Omega_0\subset\Omega.
\]
On the open arc $I$, the displacement vanishes and the exact boundary
identity gives $F(s,\xi)=R(e_\delta(\xi),s)\in M$.
Finally, $F_{\mathcal J}=s$ makes the slices distinct. Restricting
$s$ to a small complex line through $0$ gives the asserted
one-parameter family.
\end{proof}

\begin{proof}[Proof of \Cref{cl:main}]
Suppose that a nonempty locally convex $C^2$ boundary patch $M$ exists.
By \Cref{lem:no-strongly-pseudoconvex}, its Levi rank is everywhere
less than $n-1$. We divide the proof according to the maximal Levi
rank attained on $M$.
Let $r$ be the largest Levi rank attained on the relatively open
patch $M$.

If $r>0$, choose a point attaining it and a local
complex tangent frame. The Levi matrix is continuous, so a
nonzero $r\times r$ minor persists in a smaller relatively open
neighborhood. Maximality then makes the Levi rank identically
$r<n-1$ there. If $r=0$, the Levi rank is already zero throughout
$M$. Relabel the resulting nonempty open patch as $M$.

Suppose first that $r=0$.
Choose a smaller true boundary patch $M_0\subset M$ on which the
real rank of the shape operator is constant, by the same
maximal-minor argument. If its maximal real rank is zero, that
rank is already constant and no minor is needed.
By \eqref{eq:K}, the complex tangent space satisfies
$H_q=K_q\subset N_q$ for every $q\in M_0$.
By \Cref{lem:Legendre}, the constant-normal fiber through $q$
contains a neighborhood of $q$ in the real affine plane $q+N_q$.
It therefore contains a neighborhood of $q$ in the complex affine
hyperplane $q+H_q$.
Choose a complex-linear isomorphism
$L_q:\C^{n-1}\to H_q$ and $\varepsilon_q>0$ small enough that
\[
 f_q(\lambda)=q+\varepsilon_q L_q\lambda,\quad
 \lambda\in \mathbb{B}^{n-1},
\]
has image in $M_0$.
Then $f_q(0)=q$ and $\rank_\C\dd f_q(0)=n-1$.
This contradicts \Cref{thm:nonsmooth}.

We are left with $1\le r<n-1$, which already rules out $n=2$.
For $n\ge3$, apply \Cref{prop:family} to the convex local
intersection. Its nontrivial holomorphic family is attached to the
true boundary of $D$ on an open arc, contradicting
\Cref{prop:forbidden}.
\end{proof}

\smallskip

\noindent\textbf{AI Disclosure.}
All mathematical ideas and the underlying framework presented in this manuscript were developed by the authors. During the preparation and revision of the manuscript, the authors used OpenAI's ChatGPT-6.1-Sol to assist with language editing and to improve the clarity of the mathematical exposition. The authors reviewed and edited all AI-assisted content and independently verified all mathematical statements and proofs. They take full responsibility for the accuracy and integrity of the final manuscript.

\end{document}